\documentclass[10pt,a4paper]{amsart}
\usepackage{amsmath,amssymb}
\usepackage{graphicx, tikz, hyperref, subcaption, pgfplots}
\usepackage{xcolor}
\usetikzlibrary{patterns}
\usetikzlibrary{arrows.meta}
\usetikzlibrary{spy}
\pgfplotsset{compat=1.18}

\newtheorem{Proposition}{Proposition}[section]
\newtheorem{Lemma}[Proposition]{Lemma}

\newtheorem{Theorem}[Proposition]{Theorem} 
 
\newtheorem{Remark}[Proposition]{Remark}

\definecolor{TolBlue}{HTML}{005AB5}
\definecolor{TolRed}{HTML}{DC3220}
\definecolor{AGreen}{HTML}{029356}
\definecolor{ASel}{HTML}{009eb0}
\definecolor{ABlue}{HTML}{0073e6}
\definecolor{APurp}{HTML}{9b8bf4}

\newcommand{\La}{\Lambda}
\newcommand{\abs}[1]{|#1|}
\renewcommand{\leq}{\leqslant}
\renewcommand{\geq}{\geqslant}

\title[A hook and content pair bijection]{A hook and content pair bijection}

\author[S.~Albion~Ferlinc \and M.~Szwej]{Seamus Albion Ferlinc \and Micha{\l} Szwej}
\address{Faculty of Mathematics and Physics, University of Ljubljana}
\email[S.~Albion~Ferlinc]{seamus.albion@fmf.uni-lj.si}
\address{University of Bristol}
\email[M.~Szwej]{michal.szwej@bristol.ac.uk}

\date{}

\begin{document}
	
	\begin{abstract}
		Inspired by an isomorphism of $\mathrm{SL}_2(\mathbb{C})$-plethysms due to King, we give two proofs of a new hook and content pair multiset identity. Our first proof is by an explicit bijection and the second is based on the complementation symmetry for Macdonald polynomials. Along the way we give new interpretations of a bijection for hook pairs due to Bessenrodt and Krattenthaler's ``master bijection''.
	\end{abstract}
	
	\maketitle
	
	\section{Introduction}
	\subsection{Main result} Let $\lambda$ be a \emph{partition} fitting in an $r\times c$ rectangle $R$, that is, a~weakly decreasing sequence of $r$ non-negative integers $c\geq\lambda_1\geq\lambda_2\geq\cdots\geq\lambda_r\geq 0$. Also let $\lambda^\circ$ be the complement of $\lambda$ in $R$, defined by $\lambda_i^\circ=c-\lambda_{r-i+1}$ for $i=1,2,\ldots,r$. We define the \emph{Young diagram} of $\lambda$ as $[\lambda]=\{(i,j): 1\leq i\leq r, 1\leq j\leq \lambda_i\}$ and refer to its elements as \emph{boxes}. The \emph{hook pair} of $(i,j)\in[\lambda]$ is $(\lambda_i-j, \lambda_j'-i+1)$, and the \emph{content pair} is $(j-1, r-i+1)$. The sum of the elements in a hook pair gives the \emph{hook length} $h_{i,j}(\lambda)$ of $(i,j)\in[\lambda]$ and the sum of the elements in a content pair gives the \emph{$r$-shifted content}.
	Let $H(\lambda)$, $C(\lambda)$ be the multisets of hook pairs and content pairs for all boxes of its Young diagram, respectively. We construct an explicit bijection that proves the following.
	
	\begin{Theorem}\label{thm:main}
		Let $\lambda$ be a partition fitting in an $r\times c$ rectangle with complement~ $\lambda^\circ$. Then
		\[
		H(\lambda)\sqcup C(\lambda^\circ)=H(\lambda^\circ)\sqcup C(\lambda)
		\]
		as multisets.
	\end{Theorem}
	Rotating the diagram of $\lambda^\circ$ by $180^\circ$, this assertion can be conveniently visualised in two diagrams; see Figure~\ref{fig:colourful} for an example.
\begin{figure}[htbp]
	\centering
	\definecolor{pairblue}{HTML}{0072B2}
	\definecolor{pairorange}{HTML}{A34F00}

	\begin{subfigure}[t]{0.48\textwidth}
		\centering

		\begin{tikzpicture}[x=0.14\linewidth, y=0.14\linewidth,
			every node/.style={font=\large, text=black, fill=white,
				inner sep=1.4pt, outer sep=0pt}]
			\fill[pairblue!12] (0,0) rectangle (1,1);
			\path[pattern=north east lines, pattern color=pairblue]
			(0,0) rectangle (1,1);
			\fill[pairblue!12] (1,1) rectangle (2,2);
			\path[pattern=north east lines, pattern color=pairblue]
			(1,1) rectangle (2,2);
			\fill[pairblue!12] (3,2) rectangle (4,3);
			\path[pattern=north east lines, pattern color=pairblue]
			(3,2) rectangle (4,3);
			\fill[pairblue!12] (4,3) rectangle (5,4);
			\path[pattern=north east lines, pattern color=pairblue]
			(4,3) rectangle (5,4);
			\fill[pairblue!12] (6,3) rectangle (7,4);
			\path[pattern=north east lines, pattern color=pairblue]
			(6,3) rectangle (7,4);
			
			\fill[pairorange!12] (1,2) rectangle (2,3);
			\path[pattern=dots, pattern color=pairorange]
			(1,2) rectangle (2,3);
			\fill[pairorange!12] (4,2) rectangle (5,3);
			\path[pattern=dots, pattern color=pairorange]
			(4,2) rectangle (5,3);
			\fill[pairorange!12] (2,3) rectangle (3,4);
			\path[pattern=dots, pattern color=pairorange]
			(2,3) rectangle (3,4);

			\foreach \x in {1,...,6}
			{
				\draw[gray] (\x,0)--(\x,4);
			}
			\foreach \y in {1,...,3}
			{
				\draw[gray] (0,\y)--(7,\y);
			}
						\draw (0,0) rectangle (7,4);
			\draw[very thick] 		(0,0)--(1,0)--(1,1)--(2,1)--(2,2)--(4,2)--(4,3)--(5,3)--(5,4)--(7,4);
			
			\node at (0.5,0.5) {$0,1$};
			\node at (1.5,0.5) {$5,4$};
			\node at (2.5,0.5) {$4,4$};
			\node at (3.5,0.5) {$3,4$};
			\node at (4.5,0.5) {$2,4$};
			\node at (5.5,0.5) {$1,4$};
			\node at (6.5,0.5) {$0,4$};
			
			\node at (0.5,1.5) {$1,2$};
			\node at (1.5,1.5) {$0,1$};
			\node at (2.5,1.5) {$4,3$};
			\node at (3.5,1.5) {$3,3$};
			\node at (4.5,1.5) {$2,3$};
			\node at (5.5,1.5) {$1,3$};
			\node at (6.5,1.5) {$0,3$};
			
			\node at (0.5,2.5) {$3,3$};
			\node at (1.5,2.5) {$2,2$};
			\node at (2.5,2.5) {$1,1$};
			\node at (3.5,2.5) {$0,1$};
			\node at (4.5,2.5) {$2,2$};
			\node at (5.5,2.5) {$1,2$};
			\node at (6.5,2.5) {$0,2$};
			
			\node at (0.5,3.5) {$4,4$};
			\node at (1.5,3.5) {$3,3$};
			\node at (2.5,3.5) {$2,2$};
			\node at (3.5,3.5) {$1,2$};
			\node at (4.5,3.5) {$0,1$};
			\node at (5.5,3.5) {$1,1$};
			\node at (6.5,3.5) {$0,1$};
		\end{tikzpicture}
		\caption{Hook pairs in $\lambda$ and\newline
			content pairs in the complement}
	\end{subfigure}\hfill
	\begin{subfigure}[t]{0.48\textwidth}
		\centering
		\begin{tikzpicture}[x=0.14\linewidth, y=0.14\linewidth,
			every node/.style={font=\large, text=black, fill=white,
				inner sep=1.4pt, outer sep=0pt}]
			\fill[pairblue!12] (0,0) rectangle (1,1);
			\path[pattern=north east lines, pattern color=pairblue]
			(0,0) rectangle (1,1);
			\fill[pairblue!12] (1,0) rectangle (2,1);
			\path[pattern=north east lines, pattern color=pairblue]
			(1,0) rectangle (2,1);
			\fill[pairblue!12] (2,1) rectangle (3,2);
			\path[pattern=north east lines, pattern color=pairblue]
			(2,1) rectangle (3,2);
			\fill[pairblue!12] (4,2) rectangle (5,3);
			\path[pattern=north east lines, pattern color=pairblue]
			(4,2) rectangle (5,3);
			\fill[pairblue!12] (5,3) rectangle (6,4);
			\path[pattern=north east lines, pattern color=pairblue]
			(5,3) rectangle (6,4);
			
			\fill[pairorange!12] (3,0) rectangle (4,1);
			\path[pattern=dots, pattern color=pairorange]
			(3,0) rectangle (4,1);
			\fill[pairorange!12] (4,1) rectangle (5,2);
			\path[pattern=dots, pattern color=pairorange]
			(4,1) rectangle (5,2);
			\fill[pairorange!12] (6,2) rectangle (7,3);
			\path[pattern=dots, pattern color=pairorange]
			(6,2) rectangle (7,3);

			\foreach \x in {1,...,6}
			{
				\draw[gray!50] (\x,0)--(\x,4);
			}
			\foreach \y in {1,...,3}
			{
				\draw[gray!50] (0,\y)--(7,\y);
			}
			\draw (0,0) rectangle (7,4);
			\draw[very thick] 		(0,0)--(1,0)--(1,1)--(2,1)--(2,2)--(4,2)--(4,3)--(5,3)--(5,4)--(7,4);
			
			\node at (0.5,0.5) {$0,1$};
			\node at (1.5,0.5) {$0,1$};
			\node at (2.5,0.5) {$1,2$};
			\node at (3.5,0.5) {$2,2$};
			\node at (4.5,0.5) {$3,3$};
			\node at (5.5,0.5) {$4,4$};
			\node at (6.5,0.5) {$5,4$};
			
			\node at (0.5,1.5) {$0,2$};
			\node at (1.5,1.5) {$1,2$};
			\node at (2.5,1.5) {$0,1$};
			\node at (3.5,1.5) {$1,1$};
			\node at (4.5,1.5) {$2,2$};
			\node at (5.5,1.5) {$3,3$};
			\node at (6.5,1.5) {$4,3$};
			
			\node at (0.5,2.5) {$0,3$};
			\node at (1.5,2.5) {$1,3$};
			\node at (2.5,2.5) {$2,3$};
			\node at (3.5,2.5) {$3,3$};
			\node at (4.5,2.5) {$0,1$};
			\node at (5.5,2.5) {$1,2$};
			\node at (6.5,2.5) {$2,2$};
			
			\node at (0.5,3.5) {$0,4$};
			\node at (1.5,3.5) {$1,4$};
			\node at (2.5,3.5) {$2,4$};
			\node at (3.5,3.5) {$3,4$};
			\node at (4.5,3.5) {$4,4$};
			\node at (5.5,3.5) {$0,1$};
			\node at (6.5,3.5) {$1,1$};
		\end{tikzpicture}
		\caption{ Content pairs in $\lambda$ and\newline
			hook pairs in the complement}
	\end{subfigure}

	\caption{The partition $\lambda=(5,4,2,1)$ inside a $4\times 7$ rectangle. Both $(0,1)$ and $(2,2)$ appear the same number of times in the two arrays. Theorem~\ref{thm:main} asserts this for arbitrary $\lambda,R,$ and $(a,l)$.}
	\label{fig:colourful}
\end{figure}
	
	\subsection{Literature overview and the outline} Identities with multisets of hook lengths for Young diagrams arose unexpectedly in the work of Regev and Vershik on asymptotic representation theory. In \cite{RV97}, they posed several conjectures later settled by Regev and Zeilberger \cite{Regev-Zeil}, and independently Janson \cite{Janson}. 
	
	Bessenrodt generalised these results and reduced them to the following bijection, showing that hook pairs are often, perhaps surprisingly, more natural to work with than hook lengths.
	\begin{Theorem}[{\cite[Thm.~1.1]{Bessenrodt1998}}]\label{thm:bess}
	Let $\lambda$ be a partition and let $\mathcal{Q}$ be the south-east quadrant. The hook pair of $(i,j)\in[\mathcal{Q}\setminus\lambda]$ is defined as $(j-\lambda_i-1, i-\lambda_j')$, and the hook pair of $(i,j)\in[\mathcal{Q}]$ as $(j-1,i)$. Then
	\begin{equation*}
		H(\lambda)\sqcup H(\mathcal{Q})=H(\mathcal{Q}\setminus\lambda)
	\end{equation*}
	as multisets.
	\end{Theorem}
	
	In \cite{Regev-pairs}, Regev presented more hook pair multiset identities with inductive proofs, connecting them to formulas for Jack symmetric functions. They all admitted a unified proof using Krattenthaler's \emph{master bijection} \cite{Kratt}. Another bijective proof, using Dyck paths, was given by Goulden and Yong \cite{Goulden-Yong}. See also Regev's overview~\cite{Regev06} of the identities above with examples of their algebraic applications.
	
	Theorem~\ref{thm:main} extends this rich family of results.
	In fact, it may be derived as a consequence of two multiset identities of Regev.
	For a partition $\lambda$ contained in an $r\times c$ rectangle $R$ let $S_{r,c}(\lambda)$ be the skew shape $(\lambda_1+c,\dots,\lambda_r+c)/\lambda$.
	Then \cite[Corollary~1,~($a$,$\tilde{a}$)]{Regev-pairs} states that
	\begin{subequations}\label{eq:reg}
		\begin{align}\label{eq:reg1}
			H(S_{r,c}(\lambda))\sqcup C(\lambda) &= H(\lambda)\sqcup H(R), \\
			H(S_{r,c}(\lambda^\circ))\sqcup C(\lambda) &= H(\lambda)\sqcup H(R).\label{eq:reg2}
		\end{align}
	\end{subequations}
	Taking the union of the first equality with $\lambda\mapsto \lambda^\circ$ and the second we obtain
	\[
	H(\lambda)\sqcup C(\lambda^\circ)\sqcup H(S_{r,c}(\lambda^\circ))\sqcup H(R) = 
	H(\lambda^\circ)\sqcup C(\lambda)\sqcup H(S_{r,c}(\lambda^\circ))\sqcup H(R).
	\]
	Cancelling $H(S_{r,c}(\lambda^\circ))\sqcup H(R)$ yields Theorem~\ref{thm:main}. As mentioned, \eqref{eq:reg} admit bijective proofs using Krattenthaler's master bijection, however the final derivations still rely on such a cancellation as above. This problem was suggested to us by Mark Wildon, who gave inductive and algebraic proofs of the special case for hook lengths and contents in the unpublished note~\cite{Wildon}.
	
	In \S\ref{sec:bij} we give an explicit bijective proof of Theorem~\ref{thm:main}. For fixed $(a,l)$, we consider a shift $Q$ of the boundary path $P$ of $\lambda$ by $-(a+1/2, l-1/2)$. The key feature of the proof is that the signed distance function $d(t)$ between $P$ and $Q$ contains all the information about the boxes of $\lambda$ and $\lambda^\circ$ with hook pairs or content pairs equal to $(a,l)$, encoding the entire bijection. This strategy also gives a new short bijective proof of Bessenrodt's Theorem~\ref{thm:bess} and a novel interpretation of Krattenthaler's master bijection, see \S\ref{sec:master-bij}, resulting in a new, simplified perspective on the entire family of the hook and content pair multiset identities. 
	
	\subsection{\texorpdfstring{$\mathrm{SL_2}$-plethysm}{SL2-plethysm}} Our original motivation comes from the following plethystic $\mathrm{SL}_2(\mathbb{C})$-isomorphism established by King in \cite{K85}
	\begin{equation}\label{eq:iso}
		\nabla^{\lambda}\mathrm{Sym}^\ell E\cong\nabla^{\lambda^\circ}\mathrm{Sym}^\ell E,
	\end{equation}
	where $E$ is the natural representation of $\mathrm{SL}_2(\mathbb{C})$, $\nabla^\lambda$ is the Schur functor, and $\lambda^\circ$ is the complement of $\lambda$ in a rectangle with $\ell+1$ rows. It was later strengthened by Paget and Wildon in \cite{PagetWildon} by showing that the condition for the number of rows is also necessary. We will now outline the connection between this isomorphism and our main result, following \cite[\S3]{Wildon}.
	
	Equating characters at the diagonal elements $\smash{\big(\!\begin{smallmatrix} q& 0 \\ 0& q^{-1} \end{smallmatrix}\!\big)}$ of $\mathrm{SL}_2(\mathbb{C})$ in \eqref{eq:iso} yields
	\[
	s_\lambda(q^{-\ell}, q^{-\ell+2},\dots,q^{\ell})=s_{\lambda^\circ}(q^{-\ell}, q^{-\ell+2},\dots,q^{\ell}),
	\]
	where $s_\mu$ is the Schur polynomial, see \cite[\S7]{StanleyEC2}.
	Equivalently, upon replacing $q\!\mapsto\! q^{1/2}$,
	\[
	q^{-\abs{\lambda}\ell/2}s_\lambda(1,q,\dots,q^\ell)=q^{-\abs{\lambda^\circ}\ell/2}s_{\lambda^\circ}(1,q,\dots,q^\ell),
	\]
	where $\abs{\mu}=\sum_k\mu_k$. We apply Stanley's Hook Content Formula \cite[Corollary 7.21.2]{StanleyEC2} (also see \cite[\S15.3]{Stanley71}) to obtain
	\[
	q^{n(\lambda)-\abs{\lambda}\ell/2}\prod_{(i,j)\in[\lambda]}\frac{1-q^{\ell+1+j-i}}{1-q^{h_{i,j}(\lambda)}}=q^{n(\lambda^\circ)-\abs{\lambda^\circ}\ell/2}\prod_{(i,j)\in[\lambda^\circ]}\frac{1-q^{\ell+1+j-i}}{1-q^{h_{i,j}(\lambda^\circ)}},
	\]
	where $n(\mu)=\sum_{k\geq 1} (k-1)\mu_k$. Using $\lambda^\circ_k=c-\lambda_{\ell+2-k}$, where $c$ is the number of columns of the rectangle, we observe that the powers of $q$ cancel. Moreover, defining the content of $(i,j)$ to be $c_{i,j}(\lambda)=j-i$, we may write
	\[
	\prod_{(i,j)\in[\lambda]}\frac{1-q^{\ell+1+c_{i,j}(\lambda)}}{1-q^{h_{i,j}(\lambda)}}=\prod_{(i,j)\in[\lambda^\circ]}\frac{1-q^{\ell+1+c_{i,j}(\lambda^\circ)}}{1-q^{h_{i,j}(\lambda^\circ)}}.
	\]
	
	\begin{Lemma}[{\cite[Lemma~4]{Wildon}}]\label{lem:Wildon}
		Let $\mathcal{E}$ and $\mathcal{E}'$ be finite submultisets of $\mathbb{Z}_+$. In the ring $\mathbb{C}[q]$, $\prod_{e\in\mathcal{E}}(1-q^e)=\prod_{e'\in\mathcal{E}'}(1-q^{e'})$ if and only if $\mathcal{E}=\mathcal{E}'$.
	\end{Lemma}
	In light of the lemma above, after multiplying by the denominators, the last identity is a consequence of Theorem \ref{thm:main} with $r=\ell+1$ for hook lengths and $r$-shifted contents. 
	
	While this special case corresponds to the complementation symmetry of Schur polynomials, in \S\ref{sec:mac} we show that Theorem~\ref{thm:main} in full generality is a consequence of the complementation symmetry for Macdonald polynomials, combined with the $(q,t)$-analogue of the hook-content formula.
	
	\begin{Remark}
		We refer the reader to \cite{PagetWildon} for many other examples of plethystic $\mathrm{SL}_2(\mathbb{C})$-isomorphisms obtained from comparing the multisets of hook lengths and shifted contents of two partitions. We remark that \eqref{eq:iso} is the only one that generalises to hook pairs and content pairs for arbitrary parameters.
	\end{Remark}

		\section{A bijective proof of Theorem~\ref{thm:main}}\label{sec:bij}
	We say that $(i,j)\in[\lambda]$ is an $(a,l)$\emph{-hook} if its hook pair defined in the introduction equals $(a,l)$. We define an $(a,l)$\emph{-content box} for content pairs respectively.
	The reader is warned that while the first entry in a hook pair is the \emph{arm length} $a_{i,j}(\lambda):=\lambda_i-j$, the second entry is one more than the \emph{leg length} $l_{i,j}(\lambda):=\lambda_j'-i$, so that $h_{i,j}(\lambda)=a_{i,j}(\lambda)+l_{i,j}(\lambda)+1$.
	We will construct a bijection for $(a,l)$-hooks and $(a,l)$-content boxes with fixed $0\leq a<c$ and $0<l\leq r$. Note that any partition contains at most one $(a,l)$-content box.
	
\begin{figure}[h]
		\begin{tikzpicture}[scale=.6]
			
			\foreach \x in {1,...,6}
			{
				\draw[gray] (\x,0)--(\x,4);
			}
			\foreach \y in {1,...,3}
			{
				\draw[gray] (0,\y)--(7,\y);
			}
			\draw (0,0) rectangle (7,4);
			\draw[very thick, TolBlue] (0,0)--(1,0)--(1,1)--(2,1)--(2,2)--(4,2)--(4,3)--(5,3)--(5,4)--(7,4);
			\draw[very thick, dashed, TolRed, shift={(-2.5, -1.5)}] (0,0)--(1,0)--(1,1)--(2,1)--(2,2)--(4,2)--(4,3)--(5,3)--(5,4)--(7,4);
			\draw[TolRed, fill=TolRed] (6.5,2.5) circle (3pt);
			\draw[TolRed, fill=TolRed] (3.5,.5) circle (3pt);
			\draw[thick, TolBlue] (1.5,2.5) circle (3pt);
			\draw[TolRed, fill=TolRed] (4.5,1.5) circle (3pt);
			\draw[thick, TolBlue] (2.5,3.5) circle (3pt);
			\node[TolRed] at (-2,-1) {$Q$};
			\node[TolBlue] at (7.3,3.8) {$P$};
		\end{tikzpicture}
		\caption{$\lambda=(5,4,2,1)$ in a $4\times 7$ rectangle. The $(2,2)$-hooks of $\lambda$ and (rotated) $\lambda^\circ$ are marked with dots.}
		\label{fig:P-and-Q}
	\end{figure}
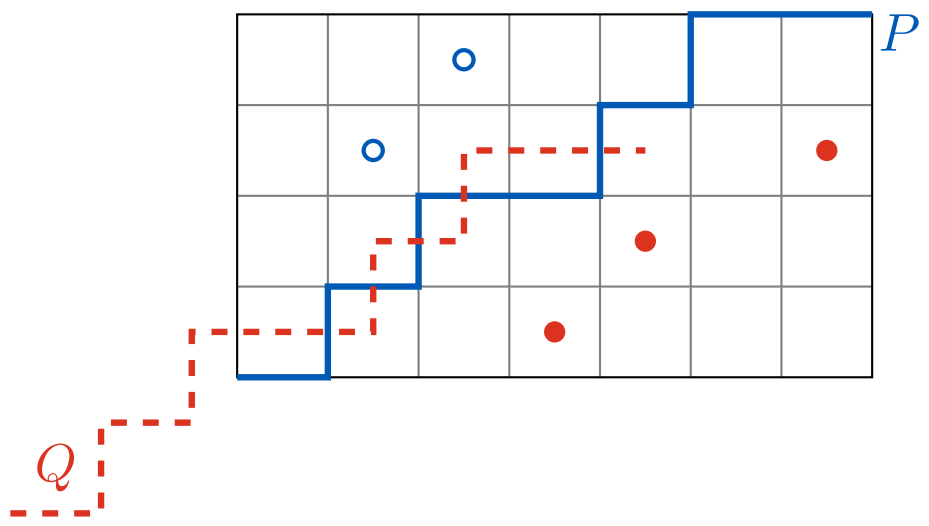

	For the rest of this section we will write $(x,y)$ for Cartesian coordinates, rather than partition-box coordinates. Let $P$ be the lattice path from $(0,0)$ to $(c,r)$ having only the boxes of $[\lambda]$ on its north-west side, see Figure~\ref{fig:P-and-Q}. Let $\alpha=(a+1/2, l-1/2)$, and $Q=P-\alpha$.
	
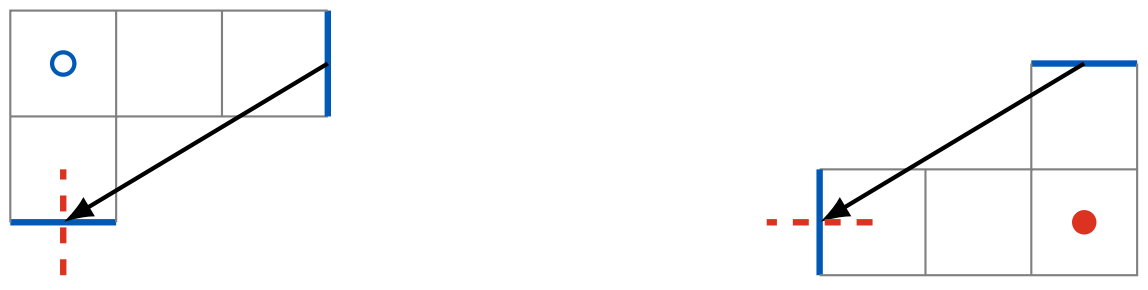
\begin{figure}[h]
		\centering
		\begin{subfigure}{.4\textwidth}
			\centering
			\begin{tikzpicture}[scale=.7]
				\draw[very thick, TolBlue] (0,0)--(1,0);
				\draw[gray] (1,0)--(1,1)--(3,1);
				\draw[very thick, TolBlue] (3,1)--(3,2);
				\draw[gray] (3,2)--(0,2)--(0,0);
				\draw[gray] (0,1)--(1,1)--(1,2);
				\draw[gray] (2,1)--(2,2);
				\draw[-{Latex[length=2mm]}, thick] (3,1.5)--(0.5,0);
				\draw[very thick, dashed, TolRed] (.5,-.5)--(.5,.5);
				\draw[thick, TolBlue] (.5,1.5) circle (3pt);
			\end{tikzpicture}
			\caption{$\lambda$-crossing}
		\end{subfigure}
		\begin{subfigure}{.4\textwidth}
			\centering
			\begin{tikzpicture}[scale=.7, rotate=180]
				\draw[very thick, TolBlue] (0,0)--(1,0);
				\draw[gray] (1,0)--(1,1)--(3,1);
				\draw[very thick, TolBlue] (3,1)--(3,2);
				\draw[gray] (3,2)--(0,2)--(0,0);
				\draw[gray] (0,1)--(1,1)--(1,2);
				\draw[gray] (2,1)--(2,2);
				\draw[{Latex[length=2mm]}-, thick] (3,1.5)--(0.5,0);
				\draw[very thick, dashed, TolRed] (2.5,1.5)--(3.5,1.5);
				\filldraw[TolRed] (.5,1.5) circle (3pt);
			\end{tikzpicture}
			\caption{$\lambda^\circ$-crossing}
		\end{subfigure}
		\caption{Hooks are in bijection with crossings}
		\label{fig:crossings}
	\end{figure}

	Notice that each crossing of $P$ and $Q$ corresponds to an $(a,l)$-hook of either $\lambda$ or $\lambda^\circ$, see Figure~\ref{fig:crossings}. When $P$ is horizontal at the crossing we call it a \emph{$\lambda$-crossing}, and when $P$ is vertical, we call it a \emph{$\lambda^\circ$-crossing}.
	
	From the diagrams above, the $(a,l)$-hooks of $\lambda$ (resp. $\lambda^\circ$) are in bijection with $\lambda$- (resp. $\lambda^\circ$-) crossings. The types of crossings alternate along the path, which will give us an alternating ordering of the $(a,l)$-hooks. Pairing the consecutive ones will later be a part of the bijection. We now make these observations more concrete.
	
\begin{figure}[h]
		\begin{tikzpicture}[scale=.75]
			\begin{axis}[
				axis x line=middle, axis y line=middle,
				ymin=-1, ymax=1.5, ytick={-1,-.5,0,.5,1}, ylabel=$d(t)$,
				xmin=0, xmax=7.5, xtick={1,...,7}, xlabel=$t$, width=\linewidth, height=4cm
				]
				\draw[very thick] (0,-.5)--(1,-.5)--(2,.5)--(3,-.5)--(4,.5)--(5,-.5)--(6,-.5)--(7,.5);
				\filldraw[TolRed] (1.5,0) circle (3pt) node[anchor=south]{\LARGE{$z_1$}};
				\draw[thick, TolBlue, fill=white] (2.5,0) circle (3pt) node[anchor=south]{\LARGE$z_2$};
				\filldraw[TolRed] (3.5,0) circle (3pt) node[anchor=south]{\LARGE$z_3$};
				\draw[thick, TolBlue, fill=white] (4.5,0) circle (3pt) node[anchor=south]{\LARGE$z_4$};
				\filldraw[TolRed] (6.5,0) circle (3pt) node[anchor=south]{\LARGE$z_5$};
			\end{axis}
		\end{tikzpicture}
		\caption{The plot of $d(t)$ with zeros $z_1<z_2<\dots<z_m$.}
		\label{fig:d-of-t}
	\end{figure}
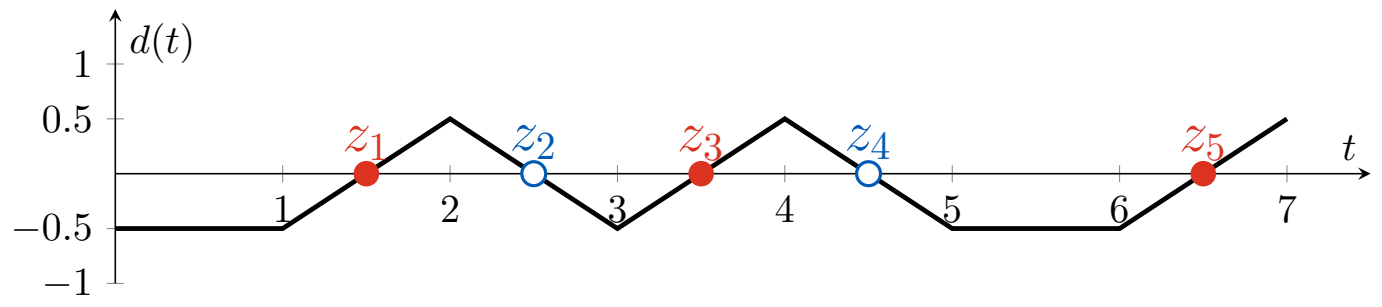

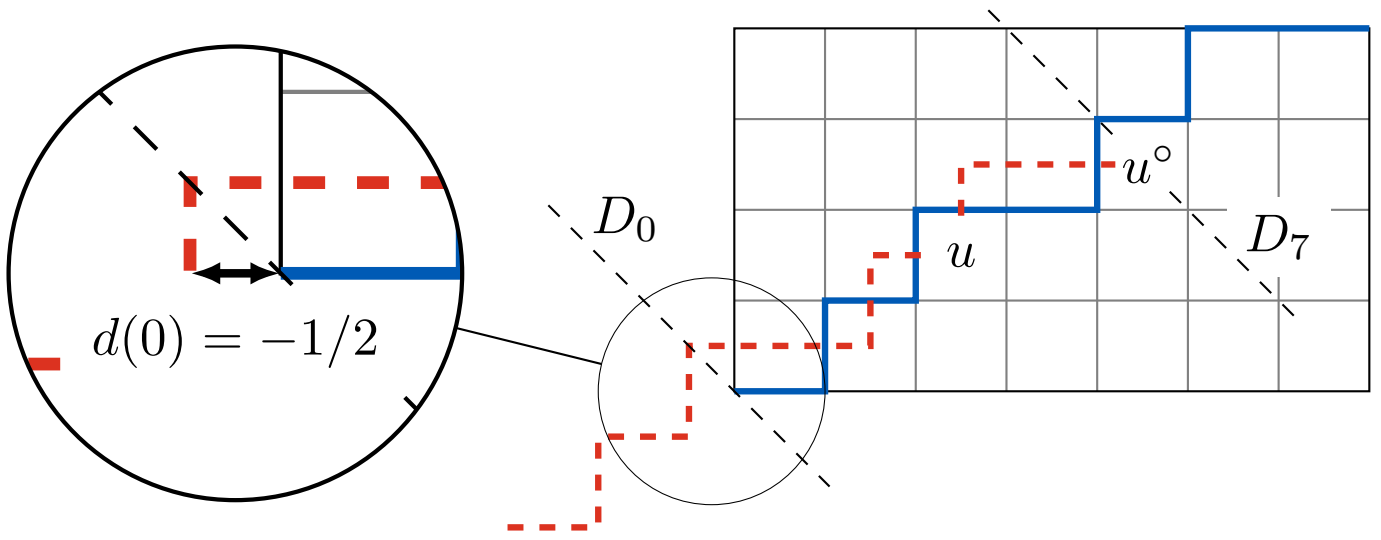
\begin{figure}[h]
		\begin{tikzpicture}[scale=.6, ]
		\begin{scope}[spy using outlines={circle, magnification=2, size=3cm, connect spies}]
				
			\foreach \x in {1,...,6}
			{
				\draw[gray] (\x,0)--(\x,4);
			}
			\foreach \y in {1,...,3}
			{
				\draw[gray] (0,\y)--(7,\y);
			}
			\draw (0,0) rectangle (7,4);
			\draw[very thick, TolBlue] (0,0)--(1,0)--(1,1)--(2,1)--(2,2)--(4,2)--(4,3)--(5,3)--(5,4)--(7,4);
			\draw[very thick, dashed, TolRed, shift={(-2.5, -1.5)}] (0,0)--(1,0)--(1,1)--(2,1)--(2,2)--(4,2)--(4,3)--(5,3)--(5,4)--(6.7,4);
			\draw[dashed] (-2.05,2.05)--(1.05,-1.05);
			\node[fill=white] at (6, 1.7) {$D_7$};
			\draw[dashed] (2.8,4.2)--(4.2,2.8);
			\draw[dashed] (4.8,2.2)--(6.2,0.8);
			
			\node at (-1.2, 1.9) {$D_0$};
			\node[fill=white] at (2.5,1.5) {$u$};
			\node at (4.57,2.5) {$u^\circ$};
			\spy on (-.15,0) in node (spy) at (-5.5,1.3);
		\end{scope}
		\node[scale=2] at (spy.center) {\begin{tikzpicture}\draw[{Latex[length=1mm]}-{Latex[length=1mm]}, thick] (-.5,0)--(-.2,0);\end{tikzpicture}};
		\node[fill=white, yshift=-4.5mm] at (spy.center) { $d(0)=-1/2$};
		\end{tikzpicture}
		\caption{The lines $D_t$, the distance measured by $d(t)$, and the boxes $u,u^\circ$ in the running example.}
		\label{fig:u}
	\end{figure}

	Consider the lines $D_t$ given by $x+y=t$ for $0\leq t\leq c+r-a-l$, see Figure~\ref{fig:u}, and define $d(t)$ as the horizontal signed distance between $P$ and $Q$ along $D_t$, concretely \[d(t)=x_{Q\cap D_t}-x_{P\cap D_t},\] where $x_S$ denotes the $x$-coordinate of point $S$, see Figure~\ref{fig:d-of-t}.	
	\vspace{-2mm}
	\begin{Lemma}\label{lem:hook-zeros}
		There is a bijection between the $(a,l)$-hooks of $\lambda$ (resp. $\lambda^\circ$) and the zeros of $d$ with negative (resp. positive) derivative.
	\end{Lemma}
	\begin{proof}
		Note that $d(t)=0$ if and only if $P$ and $Q$ intersect on $D_t$, with $d'(t)<0$ at a $\lambda$-crossing and $d'(t)>0$ at a $\lambda^\circ$-crossing. Further, by construction, there is a bijection between $\lambda$-crossings and $(a,l)$-hooks of $\lambda$, and a bijection between $\lambda^\circ$-crossings and $(a,l)$-hooks of $\lambda^\circ$. Combining the two completes the proof.
	\end{proof}
	
	Before we complete the construction, we will show that $d(t)$ also encodes the necessary information about $(a,l)$-content boxes. It is key to observe that a partition may contain at most one $(a,l)$-content box. Let $u,u^\circ\in[R]$ be the boxes of the rectangle with centres at $\alpha=(a+\frac{1}{2}, l-\frac{1}{2})$ and $(c,r)-\alpha$ respectively, see Figure~\ref{fig:u}. 
	
	It is clear that $\lambda$ contains a box with content $(a,l)$ if and only if $u\in[\lambda]$ and $\lambda^\circ$ contains an $(a,l)$-content box if and only if $u^\circ\in[R]\setminus[\lambda]$. (Recall that $[R]\setminus[\lambda]$ is the rotation of $[\lambda^\circ]$ by $180^\circ$ around the centre of $R$).
	
	The zeros of $d(t)$ lie on its domain $[0,c+r-a-l]$ and they belong to $\mathbb{Z}+\frac{1}{2}$, so both $d(0)$ and $d(c+r-a-l)$ have a (non-zero) sign. We will show that these signs encode the information about the $(a,l)$-content boxes of $\lambda$ and $\lambda^\circ$, respectively.
	
	\begin{Lemma}\label{lem:contents}
		The partition $\lambda$ contains an $(a,l)$-content box if and only if $d(0)>0$, and its complement $\lambda^\circ$ contains an $(a,l)$-content box if and only if $d(c+r-a-l)>0$.
	\end{Lemma} 
	\begin{proof}
		By definition, the shift of $-\alpha$ which sends $P$ to $Q$ also sends the centre of $u$ to $(0,0)$. The inequality $d(0)>0$ holds if and only if $P\cap D_0=(0,0)$ lies to the left of $Q\cap D_0$ or equivalently, $u$ lies to the left of $P$, that is $u\in[\lambda]$.
		
		Similarly for the second assertion: $d(c+r-a-l)>0$ if and only if the endpoint of $Q$, which is at the same time the centre of $u^\circ$, is to the right of $P\cap D_{c+r-a-l}$, that is $u^\circ\in[R]\setminus[\lambda]$. This is equivalent to $\lambda^\circ$ containing an $(a,l)$-content box.
	\end{proof}

\begin{figure}[h]
		\centering
		\begin{subfigure}{.45\textwidth}
			\centering
			\begin{tikzpicture}[scale=.6]
				\centering
				\path[pattern=north east lines, pattern color=black]
				(1,2) rectangle (2,3);
				\path[pattern=dots, pattern color=black]
				(2,3) rectangle (3,4);
				\path[pattern=vertical lines, pattern color=black]
				(4,2) rectangle (5,3);

				\foreach \x in {1,...,6}
				{
					\draw[gray] (\x,0)--(\x,4);
				}
				\foreach \y in {1,...,3}
				{
					\draw[gray] (0,\y)--(7,\y);
				}
				\draw (0,0) rectangle (7,4);
				\draw[very thick] (0,0)--(1,0)--(1,1)--(2,1)--(2,2)--(4,2)--(4,3)--(5,3)--(5,4)--(7,4);
				\node[font=\large, fill=white,
				inner sep=1.4pt, outer sep=0pt] at (1.5,2.5) {\textcolor{TolBlue}{$h_2$}};
				\node[font=\large, fill=white,
				inner sep=1.4pt, outer sep=0pt] at (2.5,3.5) {\textcolor{TolBlue}{$h_4$}};
				\node[font=\large, fill=white,
				inner sep=1.4pt, outer sep=0pt] at (4.57,2.5) {\textcolor{TolRed}{$u^\circ$}};
			\end{tikzpicture}
		\end{subfigure}
		\begin{subfigure}{.45\textwidth}
			\centering
			\begin{tikzpicture}[scale=.6]
				
				\path[pattern=north east lines, pattern color=black]
				(3,0) rectangle (4,1);
				\path[pattern=dots, pattern color=black]
				(4,1) rectangle (5,2);
				\path[pattern=vertical lines, pattern color=black]
				(6,2) rectangle (7,3);
				
				\foreach \x in {1,...,6}
				{
					\draw[gray] (\x,0)--(\x,4);
				}
				\foreach \y in {1,...,3}
				{
					\draw[gray] (0,\y)--(7,\y);
				}
				\draw (0,0) rectangle (7,4);
				\draw[very thick] (0,0)--(1,0)--(1,1)--(2,1)--(2,2)--(4,2)--(4,3)--(5,3)--(5,4)--(7,4);

				\node[font=\large, fill=white,
				inner sep=1.4pt, outer sep=0pt] at (3.5,0.5) {\textcolor{TolRed}{$h_1$}};
				\node[font=\large, fill=white,
				inner sep=1.4pt, outer sep=0pt] at (4.5,1.5) {\textcolor{TolRed}{$h_3$}};
				\node[font=\large, fill=white,
				inner sep=1.4pt, outer sep=0pt] at (6.5,2.5) {\textcolor{TolRed}{$h_5$}};
			\end{tikzpicture}
		\end{subfigure}
		\caption{The bijection above for $(2,2)$ pairs~up the boxes with the same background.}
		\label{fig:bij}
	\end{figure}
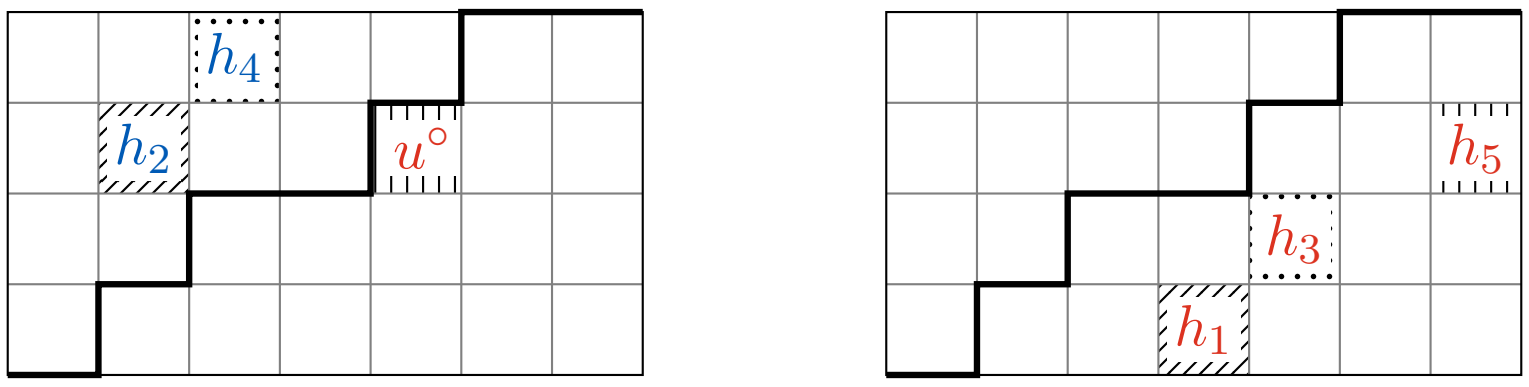

	\begin{proof}[{Proof of Theorem~\ref{thm:main}}]
		Let $z_1<\dots<z_m$ be the zeros of $d$, and let $h_i$ be the $(a,l)$-hook corresponding to $z_i$. 
		
		If $m$ is even, the sign of $d$ changes an even number of times. In particular, $d(0)$ and $d(c+r-a-l)$ are of the same sign. By Lemma~\ref{lem:contents}, either both $\lambda$ and $\lambda^\circ$ contain an $(a,l)$-content box each, in which case we pair them up in our bijection; or neither of them does. Crucially, as a consequence of Lemma~\ref{lem:hook-zeros}, the hooks $h_i$ alternate between $\lambda$ and $\lambda^\circ$, that is $h_\mathrm{even}$ belong to one partition and $h_\mathrm{odd}$ to another, see Figure~\ref{fig:bij}. Thus the pairing $(h_1, h_2),\dots, (h_{m-1},h_m)$ completes the bijection.
		
		If $m$ is odd, the sign changes an odd number of times, so by Lemma~\ref{lem:contents} precisely one of $\lambda, \lambda^\circ$ contains an $(a,l)$-content box, denoted $u$ or $u^\circ$, respectively. Thus we consider two further subcases.
		
		If $u\not\in[\lambda]$ and $u^\circ\in[R]\setminus[\lambda]$, then $d(0)<0<d(c+r-a-l)$ which means that the last zero has positive derivative, i.e. $h_m\in[\lambda^\circ]$. Thus we define our bijection as the pairing $(h_1,h_2),\ldots, (h_{m-2}, h_{m-1}), (h_m, u^\circ)$, see Figure~\ref{fig:bij}.
		
		Otherwise $u\in[\lambda]$ and $u^\circ\not\in[R]\setminus[\lambda]$, that is $d(c+r-a-l)<0<d(0)$, and so $h_m\in[\lambda]$. Similarly to the previous subcase, we define the bijection as the pairing $(h_1,h_2),\ldots,(h_{m-2},h_{m-1}),(h_m,u)$.
		
		We constructed bijections for each pair $(a,l)$ with $0\leq a <c$ and $0<l\leq r$. Since $H$ and $C$ are multiset unions of multisets of $(a,l)$-hooks and $(a,l)$-content boxes, respectively, this concludes the proof.
	\end{proof}
	
	We conclude this section by giving a short proof of Theorem~1.1 from \cite{Bessenrodt1998}.
	
	\begin{proof}[Proof of Theorem~\ref{thm:bess}]
		Let $P$ be the bi-infinite path from south to east, passing through the boundary of $\lambda$. Again, let $Q=P-\alpha$. Now $d(t)$ is defined on all of~$\mathbb{R}$ and moreover there exist $T_-, T_+$ such that $d(t)<0$ for all $t<T_-$ and $d(t)>0$ for all $t>T_+$. Thus, $d(t)$ changes sign an odd number of times, with the first and last zeros having positive derivative. By Lemma~\ref{lem:hook-zeros}, there is therefore one more $(a,l)$-hook in $[\mathcal{Q}]\setminus[\lambda]$ than in $\lambda$. The bijection sends it to the unique $(a,l)$-hook in~$\mathcal{Q}$, and pairs the remaining hooks like before.
	\end{proof}
				
	\section{Krattenthaler's master bijection}\label{sec:master-bij}
	In \cite{Kratt}, Krattenthaler shows that all of the hook and hook pair identities mentioned in the introduction, apart from Theorem~\ref{thm:main}, which is not considered there, follow from one ``master bijection''. 
	The purpose of this section is to point out that this bijection may be realised simply in terms of the ideas presented in the previous section.
	
	Let $P$ be a north-east lattice path between the lines $y=0$ and $y=r$, normalised so that the first north step occurs at $(0,0)$, and let $L(P)$ (resp.~$R(P)$) be the region to the left of $P$ (resp.~right of $P$).
	We recall the definition of hook pairs from before: label each box $s$ to the left of $P$ by $(a(s),l(s))$ where $a(s)$ is the number of boxes between $s$ and $P$ horizontally and $l(s)$ is the number of boxes between $s$ and $P$ vertically plus one. Similarly, for $s$ to the right of $P$ label it by $(a'(s),l'(s))$ where $a'(s)$ is the number of boxes between $s$ and $P$ horizontally and $l'(s)$ is the number of boxes between $s$ and $P$ vertically plus one.
	
\begin{figure}[h]
		\centering
		\begin{tikzpicture}[scale=.6, every node/.style={font=\large, text=black, fill=white,
				inner sep=1.4pt, outer sep=0pt}]
			\fill[AGreen!12] (1,2) rectangle (2,3);
			\path[pattern=north east lines, pattern color=AGreen]
			(1,2) rectangle (2,3);
			\fill[AGreen!12] (5,0) rectangle (6,1);
			\path[pattern=north east lines, pattern color=AGreen]
			(5,0) rectangle (6,1);
			
			\fill[ASel!12] (3,4) rectangle (4,5);
			\path[pattern=dots, pattern color=ASel]
			(3,4) rectangle (4,5);
			\fill[ASel!12] (6,2) rectangle (7,3);
			\path[pattern=dots, pattern color=ASel]
			(6,2) rectangle (7,3);
			
			\fill[ABlue!12] (4,5) rectangle (5,6);
			\path[pattern=vertical lines, pattern color=ABlue]
			(4,5) rectangle (5,6);
			\fill[ABlue!12] (9,5) rectangle (10,6);
			\path[pattern=vertical lines, pattern color=ABlue]
			(9,5) rectangle (10,6);

			\fill[APurp!12] (7,9) rectangle (8,10);
			\path[pattern=grid, pattern color=APurp]
			(7,9) rectangle (8,10);
			\fill[APurp!12] (10,7) rectangle (11,8);
			\path[pattern=grid, pattern color=APurp]
			(10,7) rectangle (11,8);
			
			\node at (0.5,0.5) {\textcolor{black}{$1$}};
			\node at (0.5,1.5) {\textcolor{black}{$2$}};
			\node[fill=white] at (1.5,2.5) {\textcolor{black}{$3$}};
			\node at (3.5,3.5) {\textcolor{black}{$2$}};

			\node[fill=white] at (3.5,4.5) {\textcolor{black}{$3$}};

			\node[fill=white] at (4.5,5.5) {\textcolor{black}{$3$}};
			\node at (4.5,6.5) {\textcolor{black}{$4$}};
			\node at (5.5,7.5) {\textcolor{black}{$5$}};
			\node at (7.5,8.5) {\textcolor{black}{$2$}};

			\node[fill=white] at (7.5,9.5) {\textcolor{black}{$3$}};

			\node[fill=white] at (5.5,0.5) {\textcolor{black}{$3$}};
			\node at (5.5,1.5) {\textcolor{black}{$2$}};

			\node[fill=white] at (6.5,2.5) {\textcolor{black}{$3$}};
			\node at (8.5,3.5) {\textcolor{black}{$5$}};
			\node at (8.5,4.5) {\textcolor{black}{$4$}};

			\node[fill=white] at (9.5,5.5) {\textcolor{black}{$3$}};
			\node at (9.5,6.5) {\textcolor{black}{$2$}};

			\node[fill=white] at (10.5,7.5) {\textcolor{black}{$3$}};
			\node at (12.5,8.5) {\textcolor{black}{$2$}};
			\node at (12.5,9.5) {\textcolor{black}{$1$}};
		\foreach \x in {1,...,13}
		{
			\draw[gray] (\x,0)--(\x,10);
		}
		\foreach \y in {1,...,10}
		{
			\draw[gray] (0,\y)--(13,\y);
		}
		\draw (0,0) rectangle (13,10);
		\draw[very thick] (0,0)--(3,0)--(3,2)--(4,2)--(4,3)--(6,3)--(6,5)--(7,5)--(7,7)--(8,7)--(8,8)--(10,8)--(10,10)--(13,10);
		\draw[TolRed,very thick,dashed] (1.5,-0)--(1.5,0.5)--(3.5,0.5)--(3.5,2.5)--(4.5,2.5)--(4.5,4.5)--(5.5,4.5)--(5.5,5.5)--(7.5,5.5)--(7.5,7.5)--(8.5,7.5);
		\end{tikzpicture}
		\caption{A north-east lattice path $P$ and its (truncated) shift $Q$ for $(n,m)=(2,3)$. Boxes of distance $n=2$ from $P$ have been labelled with their leg length. Shaded boxes of the same colour are mapped to one another under the bijection of Lemma~\ref{lem:master}.}
		\label{fig:Kratt}
	\end{figure}
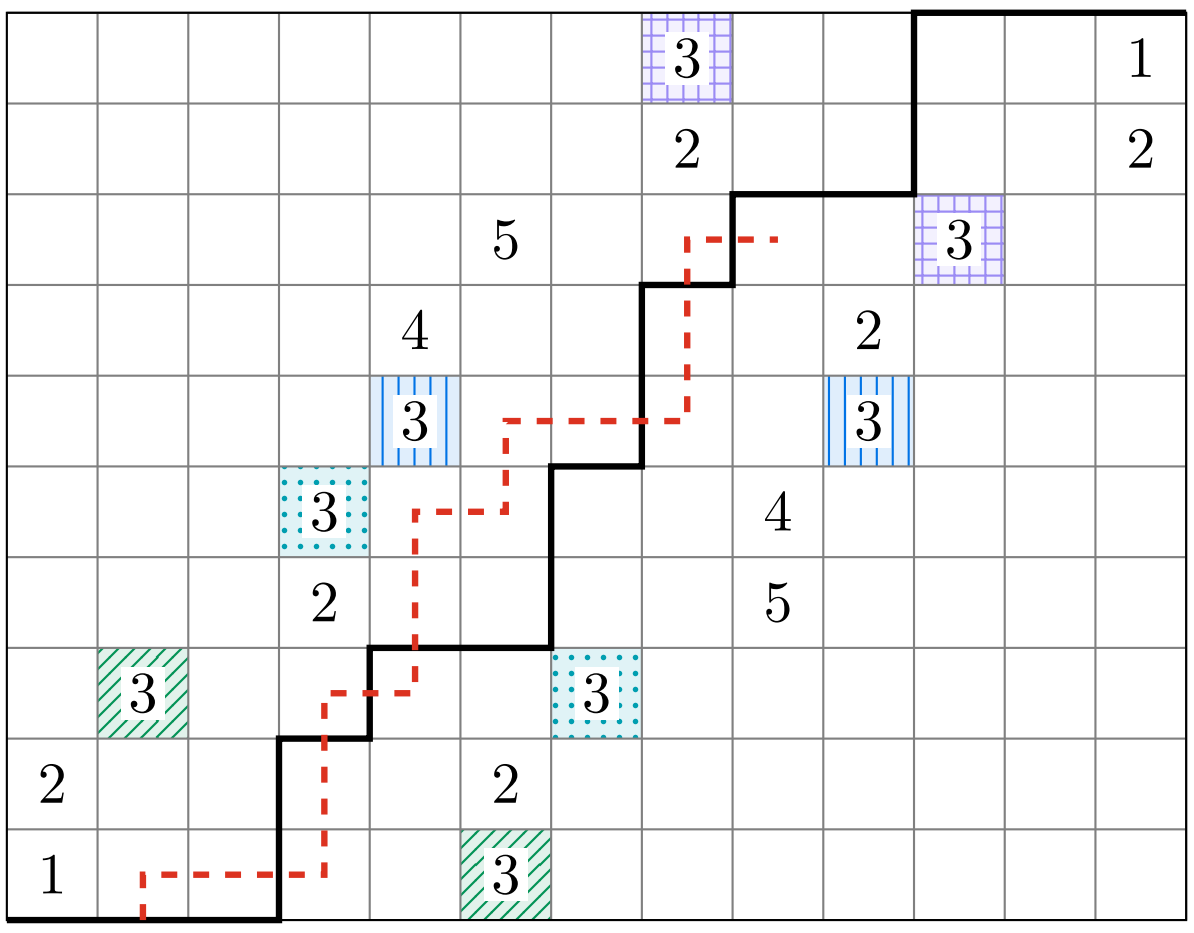

	\begin{Lemma}[{\cite[Lemma, p.~6]{Kratt}}]\label{lem:master}
		Given a lattice path $P$ as above, for each integer $n\geq 0$ there is an explicit bijection between the multisets
		\[
		\{l(s) : \text{$s\in L(P)$ and $a(s)=n$}\} \quad\text{and}\quad
		\{l'(s) : \text{$s\in R(P)$ and $a'(s)=n$}\}.
		\]
	\end{Lemma}
	\begin{proof}
		For fixed $m$ which is an element of the above sets we set $\alpha=(n+1/2,m-1/2)$ and consider the shifted lattice path $Q = P-\alpha$, where we extend both $P$ and $Q$ by horizontal steps to the left and right.
		In view of the proof of Theorem~\ref{thm:main} and Lemma~\ref{lem:hook-zeros} from the previous section, the crossings between $P$ and $Q$ alternate between $\lambda$-crossings, corresponding to $(n,m)$-hooks in $L(P)$, and $\lambda^\circ$-crossings, corresponding to $(n,m)$-hooks in $R(P)$.
		Moreover, the first crossing will be a $\lambda$-crossing and the final a $\lambda^\circ$-crossing, so that there is an even number of total crossings. Let us call the associated $(n,m)$-hooks $h_1,\dots,h_{2k}$.
		The desired bijection is then pairs the leg-lengths of the $(n,m)$-hooks $(h_{2i-1}, h_{2i})$ for $1\leq i\leq k$, running over all possible leg-lengths $m$.
		An example of the construction is given in Figure~\ref{fig:Kratt}.
	\end{proof}
	
	Krattenthaler's description of the bijection is through a sorting algorithm on increasing ``stacks'' of consecutive leg lengths. It is easy to see that our bijection and his are one and the same since both descriptions preserve the relative order of boxes with the same leg length.

	\section{A proof via the $(q,t)$-hook-content formula}\label{sec:mac}
	In this section we present a proof of the main result via the $(q,t)$-hook-content formula for Macdonald polynomials. 
	Our notations and conventions follow Macdonald's treatise \cite[\S6]{Macdonald95}.
	
	Let $\La_r$ denote the algebra of symmetric polynomials with coefficients in $\mathbb{Q}(q,t)$ in the variables $X=(x_1,\dots,x_r)$.
	The Macdonald polynomials, here denoted $P_\lambda(X;q,t)$, are a distinguished basis for $\La_r$ introduced by Macdonald \cite{Macdonald88}. 
	They reduce to the Schur polynomials for $q=t$, the Hall--Littlewood polynomials for $q=0$, and the Jack polynomials for $(q,t)=(q,q^\gamma)$ in the $q\to 1$ limit.
	Recall that on partitions of $n$ we have $\lambda>\mu$ in dominance order if
	\[
	\lambda_1+\cdots+\lambda_i\geq \mu_1+\cdots+\mu_i
	\]
	for all $1\leq i\leq \ell(\lambda)$.
	We also define the monomial symmetric polynomials as
	\[
	m_\lambda(X)=\sum_\alpha x_1^{\alpha_1}\cdots x_r^{\alpha_r},
	\]
	where the sum is over all distinct permutations $\alpha$ of the partition $\lambda$.
	It is customary to define the Macdonald polynomials $P_\lambda(X;q,t)$ to be the unique symmetric polynomials which are orthogonal with respect to the $(q,t)$-Hall inner product and are triangular in the monomial basis
	\begin{equation}\label{eq:triangularity}
		P_\lambda(X;q,t)=m_\lambda(X)+\sum_{\mu<\lambda}u_{\lambda\mu}(q,t)m_\mu(X),
	\end{equation}
	for some coefficients $u_{\lambda\mu}(q,t)\in\mathbb{Q}(q,t)$.
	For our purposes it will prove useful to use the integral inner product as defined in \cite[p.~372]{Macdonald95}.
	Let 
	\[
	\mathbb{T}^r:=\{(x_1,\dots,x_r)\in\mathbb{C}^r:\abs{x_1}=\cdots=\abs{x_r}=1\},
	\]
	be the $r$-dimensional torus. Further define for a parameter $a$ the infinite $q$-shifted factorial
	\[
	(a;q)_\infty:=\prod_{i\geq 1}(1-aq^{i-1}),
	\]
	and the density
	\[
	\Delta_r(x;q,t):=\prod_{1\leq i<j\leq r}\frac{(x_i/x_j;q)_\infty(x_j/x_i;q)_\infty}{(tx_i/x_j;q)_\infty(tx_j/x_i;q)_\infty}.
	\]
	Then for $\abs{q}, \abs{t}<1$ the inner product of interest is
	\begin{equation}\label{eq:inner}
		\langle f,g\rangle_{q,t}'
		:=\frac{1}{r!(2\pi\textup{i})^r}\int\limits_{\mathbb{T}^r}\! f(x_1,\dots,x_r)g(1/x_1,\dots,1/x_r)
		\Delta_r(x;q,t)\,\frac{\textup{d}x_1}{x_1}\cdots\frac{\textup{d}x_r}{x_r}.
	\end{equation}
	Equipped with this, we define the Macdonald polynomials $P_\lambda(X;q,t)$ to be the unique symmetric polynomials satisfying \eqref{eq:triangularity} and
	\begin{equation}\label{eq:ortho}
		\langle P_\lambda, P_\mu\rangle_{q,t}'=0\quad \text{for $\lambda\neq \mu$}.
	\end{equation}
	
	The proof of Theorem~\ref{thm:main} hinges on two standard identities for Macdonald polynomials. The first is the principal specialisation formula which, recalling $n(\mu)=\sum_{k\geq 1}(k-1)\mu_k$, is given by
	\begin{equation}\label{eq:princip}
		P_\lambda(1,t,\cdots,t^{r-1};q,t)
		=t^{n(\lambda)}\prod_{(i,j)\in[\lambda]}\frac{1-q^{j-1}t^{r-i+1}}{1-q^{a_{i,j}(\lambda)}t^{l_{i,j}(\lambda)+1}}.
	\end{equation}
	For $q=t$ this reduces to Stanley's Hook Content formula \cite[Theorem~15.3]{Stanley71}.
	The identity \eqref{eq:princip} may be found in \cite[pp.~335--337]{Macdonald95}, where a short proof based on a recursion for Pieri coefficients is provided.
	
	The second formula is the complementation symmetry; see, e.g., \cite[Equation~(4.3)]{BF99}. For $\lambda$ in an $r\times c$ rectangle it is the identity
	\begin{equation}\label{eq:compl}
		P_{\lambda}(x_1,\dots,x_r;q,t)=
		(x_1\cdots x_r)^cP_{\lambda^\circ}(1/x_1,\dots,1/x_r;q,t).
	\end{equation}
	This is a short exercise to verify using our definition of the Macdonald polynomials.
	Note that if $\lambda>\mu$ then $\lambda^{\circ}>\mu^{\circ}$. This follows from the fact that the partial sums of the complement partition are given by
	\[
	ci-\abs{\lambda}+\lambda_1+\cdots+\lambda_{r-i},
	\]
	for $1\leq i\leq r$. It is also clear that the leading term of 
	$(x_1\cdots x_r)^c P_{\lambda^\circ}(1/x_1,\dots,1/x_r;q,t)$ is $m_\lambda(X)$, and so the triangularity \eqref{eq:triangularity} holds.
	Now, substituting $(f,g)$ by
	\[
	\big((x_1\cdots x_r)^c P_{\lambda^\circ}(1/x_1,\dots,1/x_r;q,t),
	(x_1\cdots x_r)^cP_{\mu^\circ}(1/x_1,\dots,1/x_r;q,t)\big)
	\]
	in the inner product \eqref{eq:inner} leaves the product of symmetric polynomials in the integrand 
	\[
	P_{\lambda^\circ}(1/x_1,\dots,1/x_r;q,t)
	P_{\mu^\circ}(x_1,\dots,x_r;q,t).
	\]
	This establishes the orthogonality \eqref{eq:ortho} for the right-hand side of \eqref{eq:compl}, and thus the identity.
	
	\begin{proof}[Second proof of Theorem~\ref{thm:main}]
		In the complementation symmetry \eqref{eq:compl} replace $X\mapsto(1,t,\dots,t^{r-1}$), resulting in
		\[
		P_\lambda(1,t,\dots,t^{r-1};q,t)
		=t^{cr(r-1)/2}P_{\lambda^\circ}(1,t^{-1},\dots,t^{1-r};q,t).
		\]
		Since the Macdonald polynomial $P_{\lambda^\circ}$ is homogeneous of degree $\abs{\lambda^\circ}=rc-\abs{\lambda}$ we may pull out an overall factor of $t^{(1-r)(rc-\abs{\lambda})}$ and write
		\[
		P_\lambda(1,t,\dots,t^{r-1};q,t)
		=t^{cr(r-1)/2+(1-r)(rc-\abs{\lambda})}P_{\lambda^\circ}(1,t,\dots,t^{r-1};q,t).
		\]
		Applying the principal specialisation formula
		\eqref{eq:princip} then yields
		\begin{multline*}
			t^{n(\lambda)}\prod_{(i,j)\in[\lambda]}\frac{1-q^{j-1}t^{r-i+1}}{1-q^{a_{i,j}(\lambda)}t^{l_{i,j}(\lambda)+1}}\\
			=t^{n(\lambda^\circ)+cr(r-1)/2+(1-r)(rc-\abs{\lambda})}
			\prod_{(i,j)\in[\lambda^\circ]}\frac{1-q^{j-1}t^{r-i+1}}{1-q^{a_{i,j}(\lambda^\circ)}t^{l_{i,j}(\lambda^\circ)+1}}.
		\end{multline*}
		A short computation shows that the exponents of $t$ in the prefactors agree, 
		and so we have established the polynomial identity
		\begin{multline}\label{eq:multisets}
			\prod_{(i,j)\in[\lambda]}(1-q^{j-1}t^{r-i+1})
			\prod_{(i,j)\in[\lambda^\circ]}(1-q^{a_{i,j}(\lambda^\circ)}t^{l_{i,j}(\lambda^\circ)+1})\\
			=\prod_{(i,j)\in[\lambda^\circ]}(1-q^{j-1}t^{r-i+1})
			\prod_{(i,j)\in[\lambda]}(1-q^{a_{i,j}(\lambda)}t^{l_{i,j}(\lambda)+1}).
		\end{multline}
		As a final step we must show that the submultisets of exponents agree.
		To this end let $\mathcal{E}$ and $\mathcal{E}'$ be finite submultisets of $\{0,\dots,c-1\}\times\mathbb{N}$ and assume that as an identity in $\mathbb{Q}[q,t]$ we have
		\[
		\prod_{(a,l)\in\mathcal{E}}(1-q^at^l)=\prod_{(a,l)\in\mathcal{E}'}(1-q^at^l).
		\]
		The identity \eqref{eq:multisets} is precisely of this form since $a_{i,j}(\lambda)$ and $j-1$ are bounded above by $c-1$.
		Let $k$ be any integer such that $k\geq c$ and replace $(q,t)\mapsto(q,q^k)$.
		Then by Lemma~\ref{lem:Wildon}, as submultisets of $\mathbb{N}$, 
		\[
		\{a+kl : (a,l)\in\mathcal{E}\} = \{a+kl : (a,l)\in\mathcal{E}'\}.
		\]
		Now if $a+kl=a'+kl'$ then $(a-a')=k(l'-l)$, but since 
		$0\leq \abs{a-a'}\leq k-1$ we have $l=l'$ and $a=a'$. Therefore $\mathcal{E}=\mathcal{E}'$.
	\end{proof}
	
	\section*{Acknowledgements}                    We would like to thank Mark Wildon for suggesting this problem and for his support throughout, and \'Alvaro Guti\'errez for his proofreading. The first author thanks the University of Bristol for the kind hospitality during his visit in April 2026, where this project was initiated. He also thanks Matja\v{z} Konvalinka for helpful discussions. The first author was partially supported by the Austrian Science Fund (FWF) 10.55776/J5095. The second author gratefully acknowledges financial support from the Heilbronn Institute for Mathematical Research, Bristol, UK. ChatGPT was used exclusively for improving the readability of Figures~\ref{fig:colourful}, \ref{fig:bij} and \ref{fig:Kratt}.
	
	\bibliographystyle{amsplain}
	\bibliography{bibliography.bib}
\end{document}